\documentclass[UTF8]{amsart}
\usepackage{graphicx,graphics,psfrag}
\usepackage{amscd}
\usepackage{amsthm}
\usepackage{cite}
\usepackage[mathscr]{eucal}
\usepackage{amsmath,amssymb,latexsym,amsfonts}
\usepackage{soul}
\usepackage{fancyhdr}
\usepackage{float}
\usepackage[all]{xy}
\usepackage{graphicx}
\usepackage[colorlinks,linkcolor=blue, anchorcolor=blue,citecolor=blue]{hyperref}
\usepackage{easyReview}
\usepackage[ruled,vlined,linesnumbered]{algorithm2e}
\allowdisplaybreaks
\newtheorem{tpcl}{Tpcl}[section]
\newtheorem{theorem}[tpcl]{Theorem}

\theoremstyle{definition}
\newtheorem{definition}[tpcl]{Definition}
\theoremstyle{remark}
\newtheorem{remark}[tpcl]{Remark}
\newtheorem{example}[tpcl]{Example}

\numberwithin{equation}{section}

\begin{document}

\title{Gr\"{o}bner-Shirshov bases  for the Coxeter groups of the types  $G_2,F_4,E_6,E_7$ and $E_8$ }

\author[Xiaowei Pang]{Xiaowei Pang$^{*}$}
\thanks{$^{*}$ Supported by Science Foundation of Hebei Normal University No. L2022B30, Natural Science Fundation of Hebei Province (A2023205045), China Scholarship Council No.202308130193.}
\address{$^{*}$ \;Hebei Normal University\\
		Hebei Key Laboratory of Computational Mathematics and Applications\\
Hebei Research Center of the Basic Discipline Pure Mathematics}
\email{$^{*}$ \; pangxw21@hebtu.edu.cn}

\author[Jun Wang]{Jun Wang$^{**}$}
\thanks{$^{**}$ Supported by  Science and Technology Project of Hebei Education Department  No. QN2023030, Science Foundation of Hebei Normal University No. L2022B02, China Scholarship Council No.202308130195, Natural Science Fundation of Hebei Province (A2023205045).}

\address{$^{**}$ \;Hebei Normal University\\
		Hebei Key Laboratory of Computational Mathematics and Applications\\
		Hebei Center for Applied Mathematics\\
		Hebei International Joint Research Center for Mathematics and Interdisciplinary Science\\
Hebei Research Center of the Basic Discipline Pure Mathematics}
\email{$^{**}$ \; wjun@hebtu.edu.cn}


\maketitle

\begin{abstract}
 It is known that the finite Coxeter groups are classified in terms of Coxeter-Dynkin diagrams. The authors are mainly interest in the Gr\"{o}bner-Shirshov bases of finite Coxeter groups. In this paper, the authors obtain the Gr\"{o}bner-Shirshov basis for Coxeter group of type $E_8$,  while giving simpler Gr\"{o}bner-Shirshov bases for the Coxeter groups of the types $E_6,E_7$.  The new presentations of $E_n, n=6,7,8$ which are given in this paper has enormous effects in reducing the amount of calculation.
\end{abstract}

\section{Introduction}
Let $M=(m_{ij})$ be a symmetric $n\times n$ matrix over natural numbers, with $\infty$ such that $m_{ii}=1,2\leq m_{ij}\leq \infty $ for $i\ne j$. The matrix $M$ determines a group with the following presentation
$$
W=\mbox{sgp}\langle s_1,..,s_n|(s_is_j)^{m_{ij}}=1, \mbox{ where } 1\leq i,j\leq n, m_{ij}\ne \infty\rangle,
$$
which is called the Coxeter group $W=W(M)$.\par
The study of Coxeter groups can be said to go back to Greek antiquity, the symmetries and classification of regular polytopes being part of the theory. In 1935, Coxeter classified the finitie Coxeter groups in terms of Coxeter-Dynkin diagrams \cite{coxeter}. When computing in a Coxeter group, whether by hand or by machine, one faces the problem about how to efficiently represent its elements. For this purpose, we focus on  the  noncommutative Gr\"{o}bner bases theory, which is called the Gr\"{o}bner-Shirshov bases theory. The effecitve notion stems from Shirshov's composition Lemma and his algorithm \cite{MR2547481,MR0183753} for Lie algebras, and independently from Bchberger's algorithm \cite{MR2202562} of computing Gr\"{o}bner bases for commutative algebras. The Gr\"{o}bner-Shirshov bases  for some Lie algebras and superalgebras were determined in \cite{MR3899218,MR1800688,MR2028080,MR1700511,MR1414346,MR1733166,MR1432728}.  In 2006, Bokut and Vesnin \cite{MR2202556} applied Shirshov's method to some braid groups. Bokut considered the Coxeter groups of the types $A_n,B_n,D_n$ in 2001 \cite{MR1857281} and formulated a general hypothesis about Gr\"{o}bner-Shirshov bases of any Coxter group. In 2008, Lee Denis \cite{MR2440119} got the Gr\"{o}bner-Shirshov bases for the Coxeter groups of types $E_6$ and $E_7$.  In 2009, Yuqun Chen and Cihua Liu \cite{chen2009grbnershirshov} gave an counter-example for the general hypothesis by Bokut in \cite{MR1857281}. For non-crystallographic Coxeter groups, Jeong-Yup Lee and Dong-IL Lee \cite{MR3947551} calculated the Gr\"{o}bner-Shirshov bases in 2018. Thus the last piece of the Gr\"{o}bner-Shirshov bases for finite Coxeter group is the type of $E_8$,  which is more complicated than others.\par
In this paper, we rearrange the  generators of Coxeter groups to reduce calculating amount, get a new presentation of Gr\"{o}bner-Shirshov bases of the Coxeter groups of types $E_6,E_7$ and the Gr\"{o}bner-Shirshov bases of the Coxeter groups of types $E_8$. Thus the calculations of Gr\"{o}bner-Shirshov basis of Coxeter groups are beginning at 2001 by Bokut L.A. \cite{MR1857281}, up to now, Gr\"{o}bner-Shirshov bases of all  irreducible finite Coxeter group are solved. \par
The article is organized as follows. In Section \ref{pre}, we introduce the Gr\"{o}bner-Shirshov bases theory. The Gr\"{o}bner-Shirshov bases  for Coxeter groups of type $G_2,F_4,E_6,E_7,E_8$ will be set out in  Section \ref{e} and Appendix \ref{GF}.

\section{Preliminaries}\label{pre}
Let $K$ be a field,  $K\langle X\rangle$ be the non-commutative polynomial algebra generated by $X$ where $X=\{x_1,...,x_n\}$ is a finite set. Let $X^*=\{u|u=x_{i_1}\cdots x_{i_m},m\ge 0, i_1,\cdots,i_m \in \{1,\cdots,n\}\}$ be the monoid generated by $X$, and define the identity by the empty word  which is denoted by $1$.
Denote the degree of a word $w\in X^*$ by $|w|$.
\begin{definition}
A total order $\sigma$ on the set $X^*$ is called well-ordering, if every non-empty subset of $X^*$ has a minimal element under the order $\sigma$.
\end{definition}
\begin{definition}
A well-ordering $\sigma$ on the set $X^*$ is called monomial, if for every $u,v\in X^*$,
$$
u<_\sigma v \iff w_1uw_2<_\sigma w_1vw_2,\forall w_1,w_2\in X^*.
$$
\end{definition}
\begin{definition}
The deg-lex order on the $X^*$ is defined as follows: \par
$\forall ~ w_1=x_{i_1}...x_{i_m},w_2=x_{j_1}...x_{j_s}\in X^*$,
$$
w_1<_{deg-lex}w_2 \iff or\mbox{ } m<s, or \mbox{ }m=s,w_1<_{lex}w_2.
$$
\end{definition}
\begin{example}
For the set $X=\{x_1,x_2,x_3\}$, under the deg-lex order, we have
$$
1<x_1<x_2<x_3<x_1^2<x_1x_2<x_1x_3<x_2x_1<x_2^2<x_2x_3<\cdots ,
$$
obviously, under the  deg-lex order, $X^*$ is a well-ordering set, and the  deg-lex order is monomial.
\end{example}
Under the  monomial order $\sigma$ on  $X^*$, every  polynomial $f\in K\langle X\rangle$ can be  uniquely represented as
$$
f=\sum_{i}c_iw_i,
$$
where $c_i\in K\setminus\{0\},w_i\in X^*$ and  $w_i<_{\sigma}w_j,\mbox{ if } i<j$.
\begin{definition}
\begin{itemize}
\item The monomial $c_1w_1$ is called the leading term of $f$ with respect to $\sigma$ and  denoted  by $LT_\sigma(f)$. \par
\item $c_1\in  K\setminus\{0\}$ is called the leading coefficient of $f$ with respect to  $\sigma$ and denoted  by $LC_\sigma(f)$. \par
\item $w_1\in X^*$ is called the leading word of $f$ with respect to  $\sigma$ and denoted  by $LW_\sigma(f)$ or $\bar{f}$.
The polynomial $f$ is called monic if $LC_\sigma(f)=1$.
\end{itemize}
\end{definition}
\begin{definition}
For $g\in K\langle X \rangle,w\in X^*$, the $K-$module endomorphism $$r_g:K\langle X\rangle \rightarrow K\langle X\rangle $$ defined by
\begin{align*}
r_g(w)=
\begin{cases}
w-agb, &\mbox{ if }\exists~a,b\in X^*,\mbox{ s.t. }w=a\bar{g}b;\\
w,   &\mbox{ if }\forall~a,b\in X^*,\mbox{ s.t. }w\ne a\bar{g}b;
\end{cases}
\end{align*}
 is called a reduction by $g$.
\end{definition}
For a subset $S=\{g_1,\ldots,g_n\}\in K\langle X\rangle$, $f\in K\langle X\rangle$ is called a monomial which can be reduced to $0$ by $S$, if $f$ can be reduced to $0$ by a finite sequence of reduction $r_1,r_2,...,r_m(r_i\in \{r_g|g\in S\})$, then denote it by $f\xrightarrow{S}0$. If not, we denote it by $f \xrightarrow{S} \neg0$.
\begin{definition}
Under a monomial order  $\sigma$ on $X^*$, given two monic polynomials $f,g\in K\langle X \rangle $, the set $(f,g)$ defined by
\begin{align*}
(f,g)=\{fb-ag|a,b\in X^*,\bar{f}b=a\bar{g},|\bar{f}|+|\bar{g}|>|\bar{f}b|\}\cup \{f-agb|a,b\in X^*,\bar{f}=a\bar{g}b\}
\end{align*}
is called the composition set of $f$ and $g$.
\end{definition}
\begin{remark}
Generally, $(f,g)$ is not equal to $(g,f)$.
\end{remark}
Under  a monomial order $\sigma$,  let $S$ be a subset of $\{g\in K\langle X\rangle |LC(g)=1\}$, the ideal generated by $S$ is denoted by Ideal$(S)$.
\begin{definition}
$S$ is called the Gr\"{o}bner--Shirshov basis of Ideal$(S)$ under the monomial order $\sigma$, if for any $f,g\in S$, all the elements of $(f,g),(g,f)$ can  be reduced to $0$ by $S$.
\end{definition}
\begin{definition}
A Gr\"{o}bner-Shirshov basis $S=\{f_1,...,f_n\}\subset K\langle X\rangle$ under the monomial order $\sigma$ is called reduced if $\forall ~1\leq i,j\leq n,i\ne j$, $\bar{f_i}\ne a\bar{f_j}b$ for $\forall ~a,b\in X^*$.
\end{definition}
Let  $A=\mbox{sgp}\langle X|R \rangle $ be a semigroup presentation. Under a monomial order $\sigma$,
$
R=\{u_i=v_i|u_i,v_i\in X^*,v_i<_{\sigma}u_i\},
$
corresponds to a  subset
$
S=\{u_i-v_i|u_i,v_i\in X^*,v_i<_{\sigma}u_i\}\subset K\langle X\rangle.
$
The ideal generated by $S$ is denoted by Ideal$(S)$. The set $GB(A)$ got by the following Shirshov algorithm is called the Gr\"{o}bner-Shirshov basis of the semigroup $A$. \\

\begin{algorithm}[H]
\caption{Shirshov algorithm}
\KwIn{$S=\{u_i-v_i|u_i,v_i\in X^*,v_i<_{\sigma}u_i\}\subset K\langle X\rangle$}
\KwOut{GB(A)=Gr\"{o}bner-Shirshov basis of the semigroup $A$}
 $k=1, S_k=S$\;
 \While{$\exists~f,g\in S_k$, one of $(f,g), (g,f)$ can not reduced to $0$ by $S_k$}{
Let $A_k$ be the subset of  elements of $(f,g), (g,f)$ which can not reduced to $0$ by $S_k$\;
$S_{k+1}=S_{k}\cup A_k$\;
$k=k+1$\;
}
$GB(A)=S_k$.
\end{algorithm}

%
%
It's worth to mention that the number of the loop processes is not always finite and $GB(A)$ is not always a finite set.\par
%
In 2001, Bokut computed the Gr\"{o}bner-Shirshov bases of the Coxeter groups of the types $A_n,B_n,D_n$ in \cite{MR1857281}. We will set out the Gr\"{o}bner-Shirshov bases for the coxeter groups of the types $G_2,F_4,E_6,E_7$ and $E_8$ in the following sections.  It is easy to calculate the Gr\"{o}bner-Shirshov bases for the Coxeter group of type $G_2$ and $F_4$. We just list them as the  appendix without proof. \par
It is significant that  the presentation of group determines the computation of the Gr\"{o}bner-Shirshov bases under the fixed order.  In the following sections, we rearrange the generator, get a new presentation of Coxeter groups. It is worth to mention this  new presentation of group is the key step to get Gr\"{o}bner-Shirshov bases simpler and quicker than other presentations. \par

\section{$Gr\ddot{o}bner-Shirshov$ bases  of the Coxeter groups of the types  $E_6,E_7,E_8$}\label{e}
In this section, we first give the new presentations of finite Coxeter groups of types $E_6,E_7$ and $E_8$. Via this new presentations, we get the $Gr\ddot{o}bner-Shirshov$ bases  of the Coxeter groups of the types  $E_6,E_7,E_8$.
\begin{definition}\label{initial}
Define the presentations of finite Coxeter groups of types $E_n, n=6,7,8$ as follows:
\begin{equation}\label{E7-1}
\begin{aligned}
\theta_i^1&=x_i^2-1 ,&7-n\leq i\leq 5;\\
\theta_i^2&=x_{i+1}x_ix_{i+1}-x_ix_{i+1}x_i ,&7-n\leq i\leq 3;\\
\theta_{i,j}^3&=x_ix_j-x_jx_i ,&7-n\leq j+1<i\leq 4;\\
\theta_i^4&=x_{5}x_i-x_ix_{5}, &i\ne 6,5,3;\\
\theta^5&=x_{5}x_{3}x_{5}-x_{3}x_{5}x_{3};& \\
\varepsilon^1&=x_6^2-1;&\\
\varepsilon^2&=x_6x_{4}x_6-x_{4}x_6x_{4} ;&\\
\varepsilon_i^3&=x_6x_i-x_ix_6 ,&i\ne 4,6.\\
\end{aligned}
\end{equation}
\end{definition}
The Coxeter-Dynkin diagram related to this presentation  is as follows.
\begin{center}
  \begin{tikzpicture}
\draw (0.05,1) -- (0.95,1);
\draw (1.05,1) -- (1.95,1);
\draw (2.05,1) -- (2.95,1);
\draw (3.05,1) -- (3.95,1);
\draw (4.05,1) -- (4.95,1);
\draw (5.05,1) -- (5.95,1);
\draw [dotted] (6.05,1) -- (6.95,1);
\draw (7.05,1) -- (7.95,1);
\draw node at(0,1.3) {$x_6$};
\draw node at(2,-0.3) {$x_5$};
\draw node at(1,1.3) {$x_4$};
\draw node at(2,1.3) {$x_3$};
\draw node at(3,1.3) {$x_2$};
\draw node at(4,1.3) {$x_1$};
\draw node at(5,1.3) {$x_0$};
\draw node at(6,1.3) {$x_{-1}$};
\draw node at(7,1.3) {$x_{8-n}$};
\draw node at(8,1.3) {$x_{7-n}$};
\draw node at(10,1) {$n=6,7,8$};
\draw (0,1) circle(0.05);
\draw (1,1) circle(0.05);
\draw (2,1) circle(0.05);
\draw (3,1) circle(0.05);
\draw (4,1) circle(0.05);
\draw (5,1) circle(0.05);
\draw (6,1) circle(0.05);
\draw (7,1) circle(0.05);
\draw (8,1) circle(0.05);
\draw (2,0) circle(0.05);
\draw (2,0.05) -- (2,0.95);
\end{tikzpicture}
\end{center}
We denote that
\begin{align*}
x_{i,j}&=x_ix_{i-1}x_{i-2}...x_{j}, i>j; \quad x_{ii}=x_i, \quad x_{i,i+1}=1;\\
x_{i,j}^{'}&=x_ix_{i-2}x_{i-3}...x_{j},i>j+1; \quad x_{i,i-1}^{'}=x_{i},\quad x_{i,i}^{'}=1.
\end{align*}
The relations of the presentations of $E_n$ are called the initial relations of $E_n$.
\begin{theorem}\label{GBE6}
	The Gr\"{o}bner-Shirshov basis $GB(E_6)$ of the Coxeter group of the type $E_6$ is the set of the following relations together with the initial relations in Definition \ref{initial}:
	\begin{align*}
		\alpha_{i,j}^1&=x_{i+1,j}x_{i+1}-x_ix_{i+1,j} ,&1\leq j<i\leq3;\\
		\alpha_{5,j}^2&=x_{5,j}^{'}x_5-x_3x_{5,j}^{'} &1\leq j\leq 2;\\
		\beta_i^1&=x_{5,i}^{'}x_{4,i}-x_4x_{5,i}^{'}x_{4,i+1}, &1\leq i\leq3;\\
		\beta_i^2&=x_{5,i}^{'}x_4x_5-x_3x_{5,i}^{'}x_4 , &1\leq i\leq3;\\
		\beta_{i,j}^3&=x_{5,i}^{'}x_{4,j}x_{5,j}^{'}-x_3x_{5,i}^{'}x_{4,j}x_{5,j+1}^{'} ,&1\leq i<j\leq3;\\
		\eta_i&=x_{6,i}^{'}x_{5,i}^{'}-x_5x_{6,i}^{'}x_{5,i+1}^{'}, &1\leq i\leq4;\\
		\xi_{i,j}&=x_{6,i}^{'}x_{5,j}^{'}x_6-x_4x_{6,i}^{'}x_{5,j}^{'} ,&1\leq i<j\leq5;\\
		\lambda_{i,j,k}&=x_{6,i}^{'}x_{5,j}^{'}x_{4,k}x_{6,k}^{'}-x_4x_{6,i}^{'}x_{5,j}^{'}x_{4,k}x_{6,k+1}^{'} ,&1\leq i<j<k\leq4;\\
		\nu_{1,2,3,l}&=x_{6,1}^{'}x_{5,2}^{'}x_{4,3}x_5x_{6,l}^{'}x_5-x_4x_{6,1}^{'}x_{5,2}^{'}x_{4,3}x_5x_{6,l}^{'} ,&1\leq l\leq3.
	\end{align*}
\end{theorem}
\begin{proof}
The relations $\alpha_{i,j}^1,\alpha_{5,j}^2,\beta_i^1,\beta_i^2,\beta_{i,j}^3$ with the initial relations $\theta_i^1,\theta_i^2,\theta_{i,j}^3,\theta_i^4,\theta^5$ are the Gr\"{o}bner-Shirshov basis of the Coxeter group of the type $A_5$ which was proved in \cite{MR1857281}. What we need to check is the relations with the generator $x_6$.\par
Let $1\leq i^{'}\leq 4, 1\leq i<j\leq 5$, consider the composition set $(\xi_{i,j},\eta_{i^{'}})$,
\begin{align*}
(\xi_{i,j},\eta_{i^{'}})&=(x_{6,i}^{'}x_{5,j}^{'}x_6-x_4x_{6,i}^{'}x_{5,j}^{'},x_{6,i^{'}}^{'}x_{5,i^{'}}^{'}-x_5x_{6,i^{'}}^{'}x_{5,i^{'}+1}^{'})\\
&=(x_{6,i}^{'}x_{5,j}^{'}x_6-x_4x_{6,i}^{'}x_{5,j}^{'},x_6x_{4,i^{'}}x_{5,i^{'}}^{'}-x_5x_{6,i^{'}}^{'}x_{5,i^{'}+1}^{'})\\
&=(x_{6,i}^{'}x_{5,j}^{'}x_6-x_4x_{6,i}^{'}x_{5,j}^{'})x_{4,i^{'}}x_{5,i^{'}}^{'}-x_{6,i}^{'}x_{5,j}^{'}(x_6x_{4,i^{'}}^{'}x_{5,i^{'}}^{'}-x_5x_{6,i^{'}}^{'}x_{5,i^{'}+1}^{'})\\
&=x_{6,i}^{'}x_{5,j}^{'}x_5x_{6,i^{'}}^{'}x_{5,i^{'}+1}^{'}-x_4x_{6,i}^{'}x_{5,j}^{'}x_{4,i^{'}}x_{5,i^{'}}^{'}.
\end{align*}
If $j=2$, then $i=1$, let $i^{'}=1$, then
\begin{align*}
 &x_{6,i}^{'}x_{5,j}^{'}x_5x_{6,i^{'}}^{'}x_{5,i^{'}+1}^{'}-x_4x_{6,i}^{'}x_{5,j}^{'}x_{4,i^{'}}x_{5,i^{'}}^{'}\\
 =&x_{6,1}^{'}x_{5,2}^{'}x_5x_{6,1}^{'}x_{5,2}^{'}-x_4x_{6,1}^{'}x_{5,2}^{'}x_{4,1}x_{5,1}^{'}\\
 =&x_{6,1}^{'}x_3x_{5,2}^{'}x_{6,1}^{'}x_{5,2}^{'}-x_4x_{6,1}^{'}x_{5,2}^{'}x_{4,1}x_{5,1}^{'}
 (\mbox{ By the relations }\alpha_{5,j}^2)\\
 =&x_{6}x_4x_3x_2x_1x_3x_{5}x_3x_2x_{6,1}^{'}x_{5,2}^{'}-x_4x_{6}x_4x_3x_2x_1x_{5}x_3x_2x_{4,1}x_{5,1}^{'}\\
 =&x_{6}x_4x_2x_3x_2x_1x_{5}x_3x_2x_{6,1}^{'}x_{5,2}^{'}-x_4x_{6}x_4x_3x_2x_1x_{5}x_3x_2x_{4,1}x_{5,1}^{'}
 (\mbox{ By the relations }\alpha_{i,j}^1)\\
 =&x_{6}x_4x_2x_3x_2x_1x_{5}x_3x_2x_{6}x_4x_3x_2x_1x_{5}x_3x_2\\
 &-x_4x_{6}x_4x_3x_2x_1x_{5}x_3x_2x_{4}x_3x_2x_1x_{5}x_3x_2x_1\\
 =&x_{6}x_4x_2x_3x_2x_1x_{5}x_3x_2x_{6}x_4x_3x_2x_1x_{5}x_3x_2\\
 &-x_6x_{4}x_6x_3x_2x_1x_{5}x_3x_2x_{4}x_3x_2x_1x_{5}x_3x_2x_1
 (\mbox{ By the relations }\varepsilon^2)\\
 =&x_{6}x_4x_3x_2x_1x_5x_{3}x_2x_{6}x_4x_5x_3x_2x_1x_{5}x_3x_2\\
 &-x_6x_{4}x_3x_2x_1x_{5}x_3x_2x_6x_{4}x_3x_2x_1x_{5}x_3x_2x_1 (\mbox{ By the relations } \theta^2, \theta^5, \theta_i^4, \varepsilon_i^3)  \\
=&x_{6,1}^{'}x_{5,2}^{'}x_6x_4x_{5,1}^{'}x_{5,2}^{'}\\
&-x_{6,1}^{'}x_{5,2}^{'}x_6x_{4}x_{3,1}x_{5,1}^{'}.\\
\end{align*}
Actually
\begin{align*}
 &x_{5,1}^{'}x_{5,2}^{'}\\
 =&x_5x_3x_2x_1x_{5}x_3x_2\\
 =&x_5x_3x_{5}x_2x_1x_3x_2 (\mbox{ By the relations }\theta_i^4)\\
  =&x_3x_5x_{3}x_2x_1x_3x_2 (\mbox{ By the relations }\theta^5)\\
 =&x_3x_5x_2x_{3}x_2x_1x_2 (\mbox{ By the relations }\alpha_{i,j}^1)\\
 =&x_3x_5x_2x_{3}x_1x_2x_1 (\mbox{ By the relations }\theta_i^2)\\
 =&x_3x_2x_1x_5x_{3}x_2x_1 (\mbox{ By the relations }\theta_{i,j}^3)\\
 =&x_{3,1}x_{5,1}^{'}.\\
\end{align*}
The other cases for $i^{'},i,j$ are same. Thus the composition set $(\xi_{i,j},\eta_{i^{'}})$ can be reduced to $0$ by $GB(E_6)$.\par
Repeat the process, it is obtained that the composition set $(f,g), (g,f)$ for all $f,g\in GB(E_6)$ can be reduced to $0$ by $GB(E_6)$ and we complete the proof.
\end{proof}
\begin{theorem}
The Gr\"{o}bner-Shirshov basis $GB(E_7)$ of the Coxeter group of the type $E_7$ is  the set of  the following relations  together with the initial relations in Definition \ref{initial}:
\begin{align*}
\alpha_{i,j}^1&=x_{i+1,j}x_{i+1}-x_ix_{i+1,j}, &0\leq j<i\leq3;\\
\alpha_{5,j}^2&=x_{5,j}^{'}x_5-x_3x_{5,j}^{'} ,&0\leq j\leq 2;\\
\beta_i^1&=x_{5,i}^{'}x_{4,i}-x_4x_{5,i}^{'}x_{4,i+1}, &0\leq i\leq3;\\
\beta_i^2&=x_{5,i}^{'}x_4x_5-x_3x_{5,i}^{'}x_4  ,&0\leq i\leq3;\\
\beta_{i,j}^3&=x_{5,i}^{'}x_{4,j}x_{5,j}^{'}-x_3x_{5,i}^{'}x_{4,j}x_{5,j+1}^{'} ,&0\leq i<j\leq2;\\
\eta_i&=x_{6,i}^{'}x_{5,i}^{'}-x_5x_{6,i}^{'}x_{5,i+1}^{'} ,&0\leq i\leq4;\\
\xi_{i,j}&=x_{6,i}^{'}x_{5,j}^{'}x_6-x_4x_{6,i}^{'}x_{5,j}^{'} ,&0\leq i<j\leq5;\\
\lambda_{i,j,k}&=x_{6,i}^{'}x_{5,j}^{'}x_{4,k}x_{6,k}^{'}-x_4x_{6,i}^{'}x_{5,j}^{'}x_{4,k}x_{6,k+1}^{'} , &0\leq i<j<k\leq4;\\
\nu_{i,j,k,l}&=x_{6,i}^{'}x_{5,j}^{'}x_{4,k}x_5x_{6,l}^{'}x_5-x_4x_{6,i}^{'}x_{5,j}^{'}x_{4,k}x_5x_{6,l}^{'} ,&0\leq i<j<k\leq3,0\leq l\leq k;\\
\mu_{0,1,2,l,m}&=x_{6,0}^{'}x_{5,2}^{'}x_{4,2}x_{5,3}^{'}x_{6,l}^{'}x_{5,m}^{'}x_{4}-x_4x_{6,0}^{'}x_{5,1}^{'}x_{4,2}x_{5,3}^{'}x_{6,l}^{'}x_{5,m}^{'} ,&0\leq l<m\leq 2;\\
f&=x_{6,0}^{'}x_{5,1}^{'}x_{4,2}x_{5,3}^{'}x_4x_{6,i}^{'}x_{5,j}^{'}x_{4,k}x_{5,l}^{'}x_{4,m}x_{6,m}^{'} & \\
  &-x_4x_{6,0}^{'}x_{5,1}^{'}x_{4,2}x_{5,3}^{'}x_4x_{6,i}^{'}x_{5,j}^{'}x_{4,k}x_{5,l}^{'}x_{4,m}x_{6,m+1}^{'} ,&0\leq  i<j<k<l\leq 5,0\leq l<m\leq 5.
\end{align*}
\end{theorem}
\begin{proof}
  The proof is similar to the proof of Theorem \ref{GBE6}.
\end{proof}
We choose three different presentations of $E_7$, and calculate the Gr\"{o}bner-Shirshov basis by Mathematica.  One is the presentation from Definition \ref{initial}, and another two are the following (\ref{E7-2}) and (\ref{E7-3}).
\begin{equation}\label{E7-2}
\begin{aligned}
	\theta_i^1&=x_i^2-1 ,0\leq i\leq 5;\\
	\theta_i^2&=x_{i+1}x_ix_{i+1}-x_ix_{i+1}x_i ,0\leq i\leq 3;\\
	\theta_{i,j}^3&=x_ix_j-x_jx_i ,0\leq j+1<i\leq 4;\\
	\theta_i^4&=x_{5}x_i-x_ix_{5}, i\ne 6,5,3;\\
	\theta^5&=x_{5}x_{3}x_{5}-x_{3}x_{5}x_{3};\\
	\varepsilon^1&=x_6^2-1 ;\\
	\varepsilon^2&=x_6x_{4}x_6-x_{4}x_6x_{4} ;\\
	\varepsilon_i^3&=x_6x_i-x_ix_6, i\ne 4,6.\\
\end{aligned}
\end{equation}
and
\begin{equation}\label{E7-3}
\begin{aligned}
	\theta_i^1&=x_i^2-1 ,7-n\leq i\leq 5;\\
	\theta_i^2&=x_{i+1}x_ix_{i+1}-x_ix_{i+1}x_i ,7-n\leq i\leq 3;\\
	\theta_{i,j}^3&=x_ix_j-x_jx_i ,7-n\leq j+1<i\leq 4;\\
	\theta_i^4&=x_{5}x_i-x_ix_{5},i\ne 6,5,3;\\
	\theta^5&=x_{5}x_{3}x_{5}-x_{3}x_{5}x_{3};\\
	\varepsilon^1&=x_6^2-1 ;\\
	\varepsilon^2&=x_6x_{4}x_6-x_{4}x_6x_{4};\\
	\varepsilon_i^3&=x_6x_i-x_ix_6,i\ne 4,6.
\end{aligned}
\end{equation}

It is shown that for our representation (\ref{E7-1}),  it costs 120.281 seconds, another two representation (\ref{E7-2}) and (\ref{E7-3}) are 67382.8 seconds and 445.938 seconds respectively. All the computation are conducted using Mathematica on a computer with R7-4700G CPU, 8 kernels and 16G RAM.  These results illustrate that our new presentation makes the calculation of Gr\"{o}bner-Shirshov basis simpler. This new presentation is the key point to calculate the Gr\"{o}bner-Shirshov basis $GB(E_8)$, even though $E_8$ is more complicated than $E_6, E_7$. By the method  similar to the proof of Theorem \ref{GBE6}, we obtain the following Gr\"{o}bner-Shirshov basis $GB(E_8)$.

\begin{theorem}
The Gr\"{o}bner-Shirshov basis $GB(E_8)$ of the Coxeter group of the type $E_8$ is the set of the following relations  together with the initial relations in Definition \ref{initial}:
\begin{align*}
\alpha_{i,j}^1&=x_{i+1,j}x_{i+1}-x_ix_{i+1,j}, &-1\leq j<i\leq3;\\
\alpha_{5,j}^2&=x_{5,j}^{'}x_5-x_3x_{5,j}^{'} ,&-1\leq j\leq 2;\\
\beta_i^1&=x_{5,i}^{'}x_{4,i}-x_4x_{5,i}^{'}x_{4,i+1}, &-1\leq i\leq3;\\
\beta_i^2&=x_{5,i}^{'}x_4x_5-x_3x_{5,i}^{'}x_4 , &-1\leq i\leq3;\\
\beta_{i,j}^3&=x_{5,i}^{'}x_{4,j}x_{5,j}^{'}-x_3x_{5,i}^{'}x_{4,j}x_{5,j+1}^{'} ,&-1\leq i<j\leq2;\\
\eta_i&=x_{6,i}^{'}x_{5,i}^{'}-x_5x_{6,i}^{'}x_{5,i+1}^{'}, &-1\leq i\leq4;\\
\xi_{i,j}&=x_{6,i}^{'}x_{5,j}^{'}x_6-x_4x_{6,i}^{'}x_{5,j}^{'} ,&-1\leq i<j\leq5;\\
\lambda_{i,j,k}&=x_{6,i}^{'}x_{5,j}^{'}x_{4,k}x_{6,k}^{'}&\\
         &-x_4x_{6,i}^{'}x_{5,j}^{'}x_{4,k}x_{6,k+1}^{'} ,&-1\leq i<j<k\leq4;\\
\nu_{i,j,k,l}&=x_{6,i}^{'}x_{5,j}^{'}x_{4,k}x_5x_{6,l}^{'}x_5& \\
        &-x_4x_{6,i}^{'}x_{5,j}^{'}x_{4,k}x_5x_{6,l}^{'} ,&-1\leq i<j<k\leq3,-1\leq l\leq k;\\
\mu_{0,1,2,l,m}&=x_{6,0}^{'}x_{5,2}^{'}x_{4,2}x_{5,3}^{'}x_{6,l}^{'}x_{5,m}^{'}x_{4}&\\
         &-x_4x_{6,0}^{'}x_{5,1}^{'}x_{4,2}x_{5,3}^{'}x_{6,l}^{'}x_{5,m}^{'}, &-1\leq l<m\leq 2;\\
\end{align*}
\begin{align*}
\delta_{i,j,k,l,m,p}&=x_{6,i}^{'}x_{5,j}^{'}x_{4,k}x_{5,3}^{'}x_4x_{6,l}^{'}x_{5,m}^{'}x_{4,p}x_{6,k+3}^{'} &\\
          &-x_4x_{6,i}^{'}x_{5,j}^{'}x_{4,k}x_{5,3}^{'}x_4x_{6,l}^{'}x_{5,m}^{'}x_{4,p}x_{6,k+4}^{'},& \\
          &&-1\leq i<j<k\leq 2,-1\leq l<m<p\leq k+2;\\
\rho_{i,j,k,l,m}&=x_{6,i}^{'}x_{5,j}^{'}x_{4,k}x_{5,2}^{'}x_{6,l}^{'}x_{5,m}^{'}x_{4,p}x_{5,k+3}^{'} &\\
          &-x_4x_{6,i}^{'}x_{5,j}^{'}x_{4,k}x_{5,2}^{'}x_{6,l}^{'}x_{5,m}^{'}x_{4,p}x_{5,k+4}^{'},& \\
          &&-1\leq i<j<k\leq 1,-1\leq l<m<p\leq k+2;\\
\varphi_1&=x_{6,-1}^{'}x_{5,0}^{'}x_{4,1}x_{5,3}^{'}x_4x_{6,l}^{'}x_{5,m}^{'}x_{4,p}x_{5,q}^{'}x_{6,q}^{'} &\\
    &-x_4x_{6,-1}^{'}x_{5,0}^{'}x_{4,1}x_{5,3}^{'}x_4x_{6,l}^{'}x_{5,m}^{'}x_{4,p}x_{5,q}^{'}x_{6,q+1}^{'},&-1\leq  l<m<p<q\leq 4;\\
\varphi_2&=x_{6,-1}^{'}x_{5,0}^{'}x_{4,1}x_{5,2}^{'}x_4x_{6,l}^{'}x_{5,m}^{'}x_{4,p}x_{5,q}^{'}x_{6,p}^{'} &\\
    &-x_4x_{6,-1}^{'}x_{5,0}^{'}x_{4,1}x_{5,2}^{'}x_4x_{6,l}^{'}x_{5,m}^{'}x_{4,p}x_{5,q}^{'}x_{6,p+1}^{'},&-1\leq  l<m<p<q\leq 4;
\end{align*}
\begin{align*}
f_1&=x_{6,i_1}^{'}x_{5,i_2}^{'}x_{4,2}x_{5,3}^{'}x_4x_{6,j_1}^{'}x_{5,j_2}^{'}x_{4,j_3}x_{5,j_4}^{'}x_6&\\
    &-x_4x_{6,i_1}^{'}x_{5,i_2}^{'}x_{4,2}x_{5,3}^{'}x_4x_{6,j_1}^{'}x_{5,j_2}^{'}x_{4,j_3}x_{5,j_4}^{'},&-1\leq i_1<i_2\leq 1,-1\leq j_1<\cdots<j_4\leq 4;\\
f_2&=x_{6,-1}^{'}x_{5,0}^{'}x_{4,1}x_{5,2}^{'}x_{4,3}x_{6,j_1}^{'}x_{5,j_2}^{'}x_{4,j_3}x_{5,j_4}^{'}x_{6,j_2}^{'}&\\
    &-x_4x_{6,-1}^{'}x_{0}^{'}x_{4,1}x_{5,2}^{'}x_{4,3}x_{6,j_1}^{'}x_{5,j_2}^{'}x_{4,j_3}x_{5,j_4}^{'}x_{6,j_2+1}^{'},&  -1\leq j_1<\cdots<j_4\leq 4;\\
f_3&=x_{6,-1}^{'}x_{5,0}^{'}x_{4,1}x_{5,2}^{'}x_{4,3}x_5x_{6,j_1}^{'}x_{5,j_2}^{'}x_{4,j_3}x_{5,j_4}^{'}x_{6,j_1}^{'}&\\
    &-x_4x_{6,-1}^{'}x_{0}^{'}x_{4,1}x_{5,2}^{'}x_{4,3}x_5x_{6,j_1}^{'}x_{5,j_2}^{'}x_{4,j_3}x_{5,j_4}^{'}x_{6,j_1+1}^{'},&-1\leq   j_1<\cdots<j_4\leq 4;\\
f_4&=x_{6,-1}^{'}x_{5,0}^{'}x_{4,1}x_{5,2}^{'}x_{4}x_{6,j_1}^{'}x_{5,j_2}^{'}x_{4,j_3}x_{5,j_4}^{'}x_{4,j_5}x_{6,j_3}^{'}&\\
    &-x_4x_{6,-1}^{'}x_{0}^{'}x_{4,1}x_{5,2}^{'}x_{4}x_{6,j_1}^{'}x_{5,j_2}^{'}x_{4,j_3}x_{5,j_4}^{'}x_{4,j_5}x_{6,j_3+1}^{'}& -1\leq  j_1<...<j_5\leq 4\\
f_5&=x_{6,-1}^{'}x_{5,0}^{'}x_{4,1}x_{5,2}^{'}x_{4,3}x_{6,j_1}^{'}x_{5,j_2}^{'}x_{4,j_3}x_{5,j_4}^{'}x_{4,j_5}x_{6,j_2}^{'}&\\
    &-x_4x_{6,-1}^{'}x_{5,0}^{'}x_{4,1}x_{5,2}^{'}x_{4,3}x_{6,j_1}^{'}x_{5,j_2}^{'}x_{4,j_3}x_{5,j_4}^{'}x_{4,j_5}x_{6,j_2+1}^{'},& -1\leq  j_1<\cdots<j_5\leq 4;\\
f_6&=x_{6,-1}^{'}x_{5,0}^{'}x_{4,1}x_{5,3}^{'}x_{4}x_{6,j_1}^{'}x_{5,j_2}^{'}x_{4,j_3}x_{5,j_4}^{'}x_{4,j_5}x_{6,j_4}^{'}&\\
    &-x_4x_{6,-1}^{'}x_{5,0}^{'}x_{4,1}x_{5,2}^{'}x_{4,3}x_{6,j_1}^{'}x_{5,j_2}^{'}x_{4,j_3}x_{5,j_4}^{'}x_{4,j_5}x_{6,j_4+1}^{'},&-1\leq   j_1<\cdots<j_5\leq 4;\\
f_7&=x_{6,i_1}^{'}x_{5,i_2}^{'}x_{4,2}x_{5,3}^{'}x_{4}x_{6,j_1}^{'}x_{5,j_2}^{'}x_{4,j_3}x_{5,j_4}^{'}x_{4,j_5}x_{6,j_5}^{'}&\\
    &-x_4x_{6,-1}^{'}x_{5,0}^{'}x_{4,1}x_{5,2}^{'}x_{4,3}x_{6,j_1}^{'}x_{5,j_2}^{'}x_{4,j_3}x_{5,j_4}^{'}x_{4,j_5}x_{6,j_5+1}^{'},&-1\leq  i_1<i_2\leq 1, -1\leq j_1<\cdots<j_5\leq 4;\\
f_8&=x_{6,-1}^{'}x_{5,0}^{'}x_{4,1}x_{5,2}^{'}x_{4,3}x_5x_{6,j_1}^{'}x_{5,j_2}^{'}x_{4,j_3}x_{5,j_4}^{'}x_{4,j_5}x_{6,j_1}^{'}&\\
    &-x_4x_{6,-1}^{'}x_{5,0}^{'}x_{4,1}x_{5,2}^{'}x_{4,3}x_5x_{6,j_1}^{'}x_{5,j_2}^{'}x_{4,j_3}x_{5,j_4}^{'}x_{4,j_5}x_{6,j_1+1}^{'},& -1\leq  j_1<\cdots<j_5\leq 4;\\
f_9&=x_{6,i_1}^{'}x_{5,i_2}^{'}x_{4,2}x_{5,3}^{'}x_{4}x_5x_{6,-1}^{'}x_{5,0}^{'}x_{4,1}x_{5,2}^{'}x_{4,3}x_5x_{6,k_1}^{'}x_5&\\
    &-x_4x_{6,i_1}^{'}x_{5,i_2}^{'}x_{4,2}x_{5,3}^{'}x_{4}x_5x_{6,-1}^{'}x_{5,0}^{'}x_{4,1}x_{5,2}^{'}x_{4,3}x_5x_{6,k_1}^{'},& -1\leq  i_1<i_2\leq 1,-1\leq k_1\leq 3;\\
f_{10}&=x_{6,-1}^{'}x_{5,0}^{'}x_{4,1}x_{5,i_4}^{'}x_{4}x_{6,-1}^{'}x_{5,0}^{'}x_{4,1}x_{5,2}^{'}x_{4,3}x_5x_{6,k_1}^{'}x_{5,i_4}^{'}&\\
    &-x_4x_{6,i_1}^{'}x_{5,i_2}^{'}x_{4,1}x_{5,i_4}^{'}x_{4}x_{6,-1}^{'}x_{5,0}^{'}x_{4,1}x_{5,2}^{'}x_{4,3}x_5x_{6,k_1}^{'}x_{5,i_4+1}^{'},& -1\leq  1<i_4\leq 3,-1\leq k_1< i_4;\\
f_{11}&=x_{6,-1}^{'}x_{5,0}^{'}x_{4,1}x_{5,2}^{'}x_{4,3}x_{6,-1}^{'}x_{5,0}^{'}x_{4,1}x_{5,2}^{'}x_{4,3}x_5x_{6,k_1}^{'}x_{5,1}^{'}&\\
    &-x_4x_{6,-1}^{'}x_{5,0}^{'}x_{4,1}x_{5,2}^{'}x_{4,3}x_{6,-1}^{'}x_{5,0}^{'}x_{4,1}x_{5,2}^{'}x_{4,3}x_5x_{6,k_1}^{'}x_{5,2}^{'},& -1\leq  k_1\leq 0;\\
f_{12}&=x_{6,-1}^{'}x_{5,0}^{'}x_{4,1}x_{5,2}^{'}x_{4,3}x_5x_{6,-1}^{'}x_{5,0}^{'}x_{4,1}x_{5,2}^{'}x_{4,3}x_5x_{6,-1}^{'}x_{5,0}^{'}&\\
    &-x_4x_{6,-1}^{'}x_{5,0}^{'}x_{4,1}x_{5,2}^{'}x_{4,3}x_5x_{6,-1}^{'}x_{5,0}^{'}x_{4,1}x_{5,2}^{'}x_{4,3}x_5x_{6,-1}^{'}x_{5,1}^{'}.&  \\
\end{align*}
\end{theorem}

\begin{remark}
For the presentation of $E_n$ in Definition \ref{initial}, let $n=9$, it is actually the affine Coxeter group $\widetilde{E_8}$ which is an infinite group.  That means we can also  calculate the  Gr\"{o}bner-Shirshov basis of affine Coxeter group $\widetilde{E_8}$ via  the presentation in Definition \ref{initial}.
\end{remark}

\appendix
\section{ $Gr\ddot{o}bner-Shirshov$ bases  of the Coxeter groups of the type  $G_2,F_4$}\label{GF}
\begin{definition}
  Define the initial relations of $G_2$ by
\begin{equation}\label{G2}
\begin{aligned}
	\theta_i^1&=x_i^2-1, i=1,2;\\
	\theta^2&=(x_2x_1)^3-(x_1x_2)^3.
\end{aligned}
\end{equation}
\end{definition}
\begin{definition}
Define the initial relations of $F_4$ by
\begin{equation}\label{F4}
\begin{aligned}
	&\theta_i^1=x_i^2-1,i=1,2,3,4;\\
	&\theta^2=x_4x_3x_4-x_3x_4x_3;\\
	&\theta^3=x_4x_2-x_2x_4;\\
	&\theta^4=x_4x_1-x_1x_4;\\
	&\theta^5=x_3x_2x_3x_2-x_2x_3x_2x_3;\\
	&\theta^6=x_3x_1-x_1x_3;\\
	&\theta^7=x_2x_1x_2-x_1x_2x_1.
\end{aligned}
\end{equation}
\end{definition}
By the Shirshov algorithm, it is easy to obtain that
\begin{theorem}
The Gr\"{o}bner-Shirshov basis $GB(G_2)$ of the Coxeter group of the type $G_2$ is the set of the initial relations (\ref{G2}).
\end{theorem}
\begin{theorem}
The Gr\"{o}bner-Shirshov basis $GB(F_4)$ of the Coxeter group of the type $F_4$ is the set of  the  following relations together with the initial relations (\ref{F4}):
\begin{align*}
	\alpha^1&=x_{3,1}x_{3,1}-x_2x_{3,2}x_3x_1x_2;&\\
	\alpha^2&=x_{4,1}x_{3,2}x_{3,3}x_{4,j_1}x_{3,j_2}x_{3,j_3}x_{4,j_3}-x_3x_{4,2}x_{3,3}x_{4,4}x_{1,1}x_{2,2}x_{3,j_1}x_{4,j_2}x_{3,j_3}x_{4,j_4+1},
&1\leq  j_1<j_2<j_3\leq 4;\\
\alpha^3&=x_{4,1}x_{3,2}x_{3,4}x_{4,i}x_{3,3}x_{3,4}x_{4,5}-x_3x_{4,2}x_{3,3}x_{4,4}x_{1,1}x_{2,2}x_{3,i}x_{4,5}x_{3,4}x_{4,5},&i=1,2;\\
	\alpha^4&=x_{4,i}x_{3,3}x_{3,4}x_{4,3}x_{3,4}x_{3,4}x_{4,5}-x_3x_{4,2}x_{3,3}x_{4,4}x_{1,i}x_{2,3}x_{3,4}x_{4,5}x_{3,4}x_{4,5},&i=1,2;\\
	\alpha^5&=x_{4,i}x_{3,4}x_{3,4}x_{4,4}x_{3,4}x_{3,4}x_{4,5}-x_3x_{4,i}x_{3,4}x_{4,5}x_{1,2}x_{2,3}x_{3,4}x_{4,5}x_{3,4}x_{4,5},&i=1,2.\\
\end{align*}
\end{theorem}
\subsection*{Author contributions}
Jun Wang designed research; Xiaowei Pang performed research and realized algorithm;  Jun Wang and Xiaowei Pang wrote the paper.


\end{document}